\documentclass[12pt, reqno]{amsart}
\usepackage{amsmath, amsthm, amscd, amsfonts, amssymb, graphicx, color}
\usepackage[bookmarksnumbered, colorlinks, plainpages]{hyperref}
\hypersetup{colorlinks=true,linkcolor=red, anchorcolor=green,
citecolor=cyan, urlcolor=red, filecolor=magenta, pdftoolbar=true}

\newtheorem{theorem}{Theorem}[section]
\newtheorem{lemma}[theorem]{Lemma}
\newtheorem{proposition}[theorem]{Proposition}
\newtheorem{corollary}[theorem]{Corollary}
\theoremstyle{definition}

\theoremstyle{remark}
\newtheorem{remark}[theorem]{Remark}
\numberwithin{equation}{section}

\theoremstyle{theorem}
\newtheorem{other}{\bf Theorem}              

\newcommand{\hol}{{\mathcal Hol}}
\DeclareMathOperator{\og}{O} \DeclareMathOperator{\op}{o}
\DeclareMathOperator{\supp}{supp}
\newcommand{\Aut}{{\rm Aut}}
\newcommand{\h}{\mathcal{H}}

\def\D{{\mathbb D}}

\def\C{{\mathbb C}}

\newenvironment{pf}{\noindent{\emph{Proof.}}}{$\Box$ }
\newenvironment{Pf}{\noindent{\emph{Proof of}}}{$\Box$ }

\DeclareMathOperator{\defeq}{\overset{def}=}

\begin{document}

\setcounter{page}{1}

\title[A generalized Hilbert operator]{A generalized Hilbert operator acting on conformally
invariant spaces}

\author[D. Girela and N. Merch\'{a}n]{Daniel Girela$^1$$^{*}$ \MakeLowercase{and} Noel Merch\'{a}n$^1$}

\address{$^{1}$An\'alisis Matem\'atico,
Universidad de M\'alaga, Campus de Teatinos, 29071 M\'alaga, Spain.}
\email{\textcolor[rgb]{0.00,0.00,0.84}{girela@uma.es; noel@uma.es}}




\begin{abstract}
If $\mu $ is a positive Borel measure on the interval $[0, 1)$ we
let $\mathcal H_\mu $ be the Hankel matrix $\mathcal H_\mu =(\mu
_{n, k})_{n,k\ge 0}$ with entries $\mu _{n, k}=\mu _{n+k}$, where,
for $n\,=\,0, 1, 2, \dots $, $\mu_n$ denotes the moment of orden $n$
of $\mu $. This matrix induces formally the operator
$$\mathcal{H}_\mu (f)(z)=
\sum_{n=0}^{\infty}\left(\sum_{k=0}^{\infty}
\mu_{n,k}{a_k}\right)z^n$$ on the space of all analytic functions
$f(z)=\sum_{k=0}^\infty a_kz^k$, in the unit disc $\D $. This is a
natural generalization of the classical Hilbert operator. The action
of the operators $H_{\mu }$ on Hardy spaces has been recently
studied. This paper is devoted to study the operators $H_\mu $
acting on certain conformally invariant spaces of analytic functions
on the disc such as the Bloch space, $BMOA$, the analytic Besov
spaces, and the $Q_s$ spaces.
\end{abstract}

\maketitle
\section{Introduction}\label{intro}

Let $\D=\{z\in\C: |z|<1\}$ denote the open unit disc in the complex
plane $\C$ and let $\hol (\D)$ be the space of all analytic
functions in $\D$ endowed with the topology of uniform convergence
in compact subsets. We also let $H^p$ ($0<p\le \infty $) be the
classical Hardy spaces. We refer to \cite{D} for the notation and
results regarding Hardy spaces.
\par
If\, $\mu $ is a finite positive Borel measure on $[0, 1)$ and $n\,
= 0, 1, 2, \dots $, we let $\mu_n$ denote the moment of order $n$ of
$\mu $, that is, $\mu _n=\int _{[0,1)}t^n\,d\mu (t),$ and we let
$\mathcal H_\mu $ be the Hankel matrix $(\mu _{n,k})_{n,k\ge 0}$
with entries $\mu _{n,k}=\mu_{n+k}$. The matrix $\mathcal H_\mu $
induces formally an operator, which will be also called $\mathcal
H_\mu $, on spaces of analytic functions by its action on the Taylor
coefficients: \,$ a_n\mapsto \sum_{k=0}^{\infty} \mu_{n,k}{a_k},
\quad n=0,1,2, \dots . $\newline To be precise, if\,
 $f(z)=\sum_{k=0}^\infty a_kz^k\in \hol (\D )$
we define
\begin{equation*}\label{H}
\mathcal{H}_\mu (f)(z)= \sum_{n=0}^{\infty}\left(\sum_{k=0}^{\infty}
\mu_{n,k}{a_k}\right)z^n,
\end{equation*}
whenever the right hand side makes sense and defines an analytic
function in $\D $.
\par\medskip If $\mu $ is the Lebesgue measure on $[0,1)$ the matrix
$\mathcal H_\mu $ reduces to the classical Hilbert matrix \,
$\mathcal H= \left ({(n+k+1)^{-1}}\right )_{n,k\ge 0}$, which
induces the classical Hilbert operator $\h$ which has extensively
studied recently (see \cite{AlMonSa,Di,DiS,DJV,LNP}).
\par Galanopoulos and Pel\'{a}ez \cite{Ga-Pe2010} described
the measures $\mu $ so that the generalized Hilbert operator
$\mathcal H_\mu $ becomes well defined and bounded on $H^1$.
Chatzifountas, Girela and Pel{\'a}ez \cite{Ch-Gi-Pe} extended this
work describing those measures $\mu$ for which $\h_\mu $ is a
bounded operator from $H^p$ into $H^q$,  $0<p,q<\infty $.
\par
Obtaining an integral representation of $H_\mu$ plays a basic role
in these works. If $\mu $ is as above, we shall write throughout the
paper
\begin{equation}\label{Imu}I_\mu
(f)(z)=\int_{[0,1)}\frac{f(t)}{1-tz}\,d\mu (t),\end{equation}
whenever the right hand side makes sense and defines an analytic
function in $\D $. It turns out that the operators $H_\mu $ and
$I_\mu$ are closely related. In fact, in \cite{Ga-Pe2010} and
\cite{Ch-Gi-Pe} the measures
 $\mu $ for which the operator $I_\mu $ is
well defined in $H^p$ ($0<p<\infty $) are characterized and it is
proved that for such measures we have $\h_\mu (f)=I_\mu (f)$ for all
$f\in H^p$. These measures are Carleson-type measures.
\par If $I\subset \partial\D$ is an
arc, $\vert I\vert $ will denote the length of $I$. The
\emph{Carleson square} $S(I)$ is defined as
$S(I)=\{re^{it}:\,e^{it}\in I,\quad 1-\frac{|I|}{2\pi }\le r <1\}$.
\par If $\, s>0$ and $\mu$ is a positive Borel  measure on  $\D$,
we shall say that $\mu $
 is an $s$-Carleson measure
  if there exists a positive constant $C$ such that
\[
\mu\left(S(I)\right )\le C{|I|^s}, \quad\hbox{for any interval
$I\subset\partial\D $}.
\]
\par
If $\mu $ satisfies $\displaystyle{\lim_{\vert I\vert\to 0}\frac{\mu
\left (S(I)\right )}{\vert I\vert ^s}=0}$, then we say that $\mu $
is a\, {\it vanishing $s$-Carleson measure}. \par  A $1$-Carleson
measure, respectively, a vanishing $1$-Carleson measure, will be
simply called a Carleson measure, respectively, a vanishing Carleson
measure.
\par\medskip We recall that Carleson \cite{Carl}
proved that $H^p\,\subset \,L^p(d\mu )$ ($0<p<\infty $), if and only
if $\mu $ is a Carleson measure. This result was extended by Duren
\cite{Du:Ca} (see also \cite[Theorem\,\@9.\,\@4]{D}) who proved that
for $0<p\le q<\infty $, $H^p\subset L^q(d\mu )$ if and only if $\mu
$ is a $q/p$-Carleson measure.
\par\medskip
Following \cite{Zhao}, if $\mu$ is a positive Borel measure on $\D$,
$0\le \alpha <\infty $, and $0<s<\infty $ we say that $\mu$ is an
 $\alpha$-logarithmic $s$-Carleson measure if there exists a positive
 constant $C$ such that
 \[\frac{
\mu\left(S(I)\right )\left(\log \frac{2\pi }{\vert I\vert }\right
)^\alpha }{|I|^s}\le C, \quad\hbox{for any interval
$I\subset\partial\D $}.
\]
If $\mu\left(S(I)\right )\left(\log \frac{2\pi}{\vert I\vert }\right
)^\alpha \,=\,\op \left (|I|^s\right )$, as $\vert I\vert \to 0$, we
say that $\mu $ is a vanishing $\alpha$-logarithmic $s$-Carleson
measure.

\par\medskip
A positive Borel measure $\mu $ on $[0, 1)$ can be seen as a Borel
measure on $\mathbb D$ by identifying it with the measure $\tilde
\mu $ defined by $$ \tilde \mu (A)\,=\,\mu \left (A\cap [0,1)\right
),\quad \text{for any Borel subset $A$ of $\mathbb D$}.$$  In this
way a positive Borel measure $\mu $ on $[0, 1)$ is an $s$-Carleson
measure if and only if there exists a positive constant $C$ such
that
\[
\mu\left([t,1)\right )\le C(1-t)^s, \quad 0\le t<1,
\]
and we have similar statements for vanishing $s$-Carleson measures
and for $\alpha$-logarithmic $s$-Carleson and vanishing
$\alpha$-logarithmic $s$-Carleson measures.

\par\medskip Our main
aim in this paper is studying the operators $\mathcal H_\mu $ acting
on conformally invariant spaces.
\par
It is a standard fact that the set of all disc automorphisms ({\it
i.e.\/}, of all one-to-one analytic maps $f $ of $\mathbb D $ onto
itself), denoted $\Aut(\mathbb D) $, coincides with the set of all
M\"obius transformations of $\mathbb D $ onto itself:
$$ \Aut(\mathbb D  )=\{ \lambda \varphi_a : \vert a \vert
<1, \vert \lambda \vert =1 \} \,, $$ where $\varphi_
a(z)=(a-z)/(1-\overline a z) $.
\par
A space $X $ of analytic functions in $\D $, defined via a semi-norm
$\rho $, is said to be {\it conformally invariant\/} or {\it
M\"obius invariant\/} if whenever $f\in X $, then also $f\circ
\varphi \in X $ for any $\varphi\in\operatorname{Aut}(\D) $ and,
moreover, $\rho (f\circ\varphi )\le C\rho (f) $ for some positive
constant $C $ and all $f\in X $. A great deal of information on
conformally invariant spaces can be found in \cite{AFP, DGV1, RT}.
\par Let us start considering the Bloch space and $BMOA$.
The {\it Bloch space\/} $\mathcal B $ consists of all analytic
functions $f $ in $\mathbb D $ with bounded invariant derivative:
$$f\in \mathcal B \,\,\,\Leftrightarrow\,\,\,\Vert f\Vert _
{\mathcal B }\defeq \vert f(0)\vert +\sup\sb {z \in {\mathbb D
}}\,(1-|z|\sp2)\,|f\sp\prime(z)|<\infty \,. $$ The {\it little Bloch
space\/} $\mathcal B_0 $ is the closure of the polynomials in the
above norm of $\mathcal B $ and consists of all functions $f $
analytic in $\mathbb D $ for which $$ \lim_ {|z|\to 1} (1-|z|\sp2)
|f\sp\prime(z)| =0 \,.
$$ A classical source for the Bloch space is \cite{ACP}; see also
\cite{Zhu-book}. Rubel and Timoney \cite{RT} proved that $\mathcal
B$ is the biggest \lq\lq natural\rq\rq\, conformally invariant
space.
\par The space $BMOA $ consists of those functions $f $ in $H\sp1 $
whose boundary values have bounded mean oscillation on the unit
circle $\partial\mathbb D $ as defined by F.~John and L.~Nirenberg.
There are many characterizations of $BMOA $ functions. Let us
mention the following:
\par\smallskip
{\it If $f $ is an analytic function in $\mathbb D$, then $f\in BMOA
$ if and only if $$ \Vert f\Vert _ {BMOA}\defeq\vert f(0)\vert
+\,\Vert f\Vert_{\star }\,<\,\infty ,$$ where
$$\Vert f\Vert_{\star }
\,\defeq \,\sup_ {a\in \mathbb D }\Vert f\circ \varphi \sb
a-f(a)\Vert _ {H\sp 2}. $$} It is clear that the seminorm $\Vert
\cdot\Vert_{\star }$ is conformally invariant. If
$$\lim_{\vert a\vert \to 1}\Vert f\circ
\varphi _ a-f(a)\Vert _ {H\sp 2}=0$$ we say that $f$ belongs to the
space $VMOA$. We mention \cite{Ba, G:BMOA} as general references for
the spaces $BMOA$ and $VMOA$. Let us recall that $$H^\infty
\subsetneq BMOA \subsetneq  \bigcap_{0<p<\infty} H^p\quad
\text{and}\,\, BMOA\subsetneq\mathcal B.$$
\par\medskip Other important M\"{o}bius invariant spaces are the
analytic Besov spaces $B^p$ ($1< p<\infty $) and the $Q_s$-spaces
$(s>0)$. These spaces will be considered in
Section\,\@\ref{section-Besov-Qs}.
\par\medskip
 We close this section
noticing that, as usual, we shall be using the convention that
$C=C(p, \alpha ,q,\beta , \dots )$ will denote a positive constant
which depends only upon the displayed parameters $p, \alpha , q,
\beta \dots $ (which sometimes will be omitted) but not  necessarily
the same at different occurrences. Moreover, for two real-valued
functions $E_1, E_2$ we write $E_1\lesssim E_2$, or $E_1\gtrsim
E_2$, if there exists a positive constant $C$ independent of the
arguments such that $E_1\leq C E_2$, respectively $E_1\ge C E_2$. If
we have $E_1\lesssim E_2$ and  $E_1\gtrsim E_2$ simultaneously then
we say that $E_1$ and $E_2$ are equivalent and we write $E_1\asymp
E_2$.
\bigskip
\section{The operator $\mathcal H_\mu $ acting on $BMOA$ and the Bloch
space} We start characterizing those $\mu $ for which the operator
$I_\mu $ is well defined in $BMOA$ and in the Bloch space. It turns
out that they coincide. \begin{theorem}\label{ImuBBMOA} Let $\mu $
be a positive Borel measure on $[0,1)$. Then the following
conditions are equivalent:
\begin{itemize}
\item[(i)] $\int_{[0,1)}\log \frac{2}{1-t}\,d\mu (t)\,<\,\infty $.
\item[(ii)] For any given $f\in \mathcal B$, the integral in (\ref{Imu})
converges for all $z\in \mathbb D$ and the resulting function $I_\mu
(f)$ is analytic in \,$\D $.
\item[(iii)] For any given $f\in BMOA$, the integral in (\ref{Imu})
converges for all $z\in \mathbb D$ and the resulting function $I_\mu
(f)$ is analytic in \,$\D $.
\end{itemize}
\end{theorem}
\begin{pf}\par (i)\,$\Rightarrow $\, (ii). It is well known that there exists a positive
constant $C$ such that
\begin{equation}\label{gr-Bloch}\vert f(z)\vert \le C\Vert f\Vert _{\mathcal
B}\,\log \frac{2}{1-\vert z\vert },\quad (z\in \mathbb
D),\,\,\,\text{for every $f\in \mathcal B$},\end{equation} (see
\cite[p.\,\@13]{ACP}). Assume (i) and set $A=\int_{[0,1)}\log
\frac{2}{1-t}\,d\mu (t)$. Using (\ref{gr-Bloch}) we see that
\begin{equation}\label{int-f-B}\int_{[0,1)}\vert f(t)\vert \,d\mu
(t)\,\le C\Vert f\Vert_{\mathcal B}
\int_{[0,1)}\log\frac{2}{1-t}\,d\mu (t)\,=\,AC\Vert f\Vert_{\mathcal
B},\quad f\in \mathcal B.\end{equation} This implies that
\begin{equation}\label{Imuabs}\int_{[0,1)}\frac{\vert f(t)\vert
}{\vert 1-tz\vert }\,d\mu (t)\le \frac{AC\Vert f\Vert_{\mathcal
B}}{1-\vert z\vert },\quad (z\in \mathbb D),\quad f\in \mathcal
B.\end{equation} Using (\ref{int-f-B}), (\ref{Imuabs}), and Fubini's
theorem we see that if $f\in \mathcal B$ then:
\begin{itemize}\item For every $n\in \mathbb N$, the integral $\int_{[0,1)}t^nf(t)\,d\mu (t)$ converges
absolutely and $$\sup_{n\ge 0}\left\vert \int_{[0,1)}t^nf(t)\,d\mu
(t)\right \vert <\infty .$$
\item The integral $\int_{[0,1)}\frac{f(t)
}{1-tz}\,d\mu (t)$ converges absolutely, and
\begin{equation*}\int_{[0,1)}\frac{f(t) }{1-tz}\,d\mu (t)
\,=\,\sum_{n=0}^\infty \left (\int_{[0,1)}\,t^nf(t)\,d\mu
  (t)\right )z^n,\quad z\in \D .
  \end{equation*} \end{itemize}
  Thus, if $f\in \mathcal B$ then $I_\mu (f)$ is a well defined analytic function in $\mathbb D
  $ and \begin{equation*}I_\mu (f)(z)\,=\,\sum_{n=0}^\infty \left (\int_{[0,1)}\,t^nf(t)\,d\mu
  (t)\right )z^n,\quad z\in \D .\end{equation*}
  \par\medskip The implication (ii)\, $\Rightarrow$ \, (iii) is
  clear because $BMOA\subset \mathcal B$.
  \par\medskip (iii)\,$\Rightarrow $\, (i).
  Suppose (iii). Since the function $F(z)=\log \frac{2}{1-z}$
  belongs to $BMOA$, $I_\mu (F)(z)$ is well defined for every $z\in
  \mathbb D$. In particular $$I_\mu
  (F)(0)=\int_{[0,1)}\log\frac{2}{1-t}\,d\mu (t)$$ is a complex
  number. Since $\mu $ is a positive measure and $\log\frac{2}{1-t}>0$
  for all $t\in [0,1)$, (i) follows.
 \end{pf}
\par\bigskip Our next aim is characterizing the measures $\mu $ so
that $I_\mu $ is bounded in $BMOA$ or $\mathcal B$ and seeing
whether or not $I_\mu $ and $H_\mu $ coincide for such measures. We
have the following results.
\par\bigskip
\begin{theorem}\label{ImuboundedBMOA} Let $\mu $ be a positive Borel measure on $[0,1)$
with $\int_{[0,1)}\log \frac{2}{1-t}d\mu (t)<\,\infty $. Then the
following three conditions are equivalent:
\begin{itemize}\item[(i)] The measure $\nu $ defined by $d\nu
(t)=\log \frac{2}{1-t}\,d\mu (t)$ is a Carleson measure.
\item[(ii)] The operator $I_\mu $ is bounded from $\mathcal B$ into
$BMOA$.
\item[(iii)] The operator $I_\mu $ is bounded from $BMOA$ into
itself.
\end{itemize}
\end{theorem}
\par\bigskip
\begin{theorem}\label{Imu=Hmu}
Let $\mu $ be a positive Borel measure on $[0,1)$ with
$\int_{[0,1)}\log \frac{2}{1-t}\,d\mu (t)\,<\,\infty $. If the
measure $\nu $ defined by $d\nu (t)=\log \frac{2}{1-t}\,d\mu (t)$ is
a Carleson measure, then $\mathcal H_\mu $ is well defined on the
Bloch space and
\begin{equation*}\label{HmuImu} \mathcal H_\mu (f)\,=\,I_\mu (f),\quad
\text{for all $f\in \mathcal B$}.\end{equation*}
\end{theorem}

\par\bigskip
Theorem\,\@\ref{ImuboundedBMOA} and Theorem\,\@\ref{Imu=Hmu}
together yield the following.
\begin{theorem}\label{HmuboundedBMOA} Let $\mu $ be a positive Borel measure on
$[0,1)$ such that the measure $\nu $ defined by $d\nu (t)=\log
\frac{2}{1-t}\,d\mu (t)$ is a Carleson measure. Then the operator
$\mathcal H_\mu $ is bounded from $\mathcal B$ into $BMOA$.
\end{theorem}
\par\bigskip
\begin{Pf}{\,\em{Theorem \ref{ImuboundedBMOA}.}}
Since $\int_{[0,1)} \log\frac{2}{1-t} d\mu(t)<\infty$,
(\ref{gr-Bloch}) implies that
$$\int_{[0,1)}\vert f(t)\vert\,d\mu (t)<\infty ,\quad \text{for all
$f\in \mathcal B$}$$ and this implies that
\begin{equation*}
 \int_0^{2\pi}\int_{[0,1)} \left|\frac{f(t)g(e^{i\theta})}{1-re^{i\theta}t}\right| d\mu(t)d\theta<\infty,
 \quad 0\le r<1,\,f\in\mathcal B,\,g\in H^1.
\end{equation*}
Using this, Fubini's theorem and Cauchy's integral representation of
$H^1$-functions \cite[Theorem\,\@3.\,\@6]{D}, we deduce that
whenever $f\in\mathcal B$ and $g\in H^1$ we have
\begin{align}\label{duality1}
& \int_0^{2\pi}
I_\mu(f)(re^{i\theta})\overline{g(e^{i\theta})}\,d\theta
\,=\,\int_0^{2\pi} \left(\int_{[0,1)}\frac{f(t)
d\mu(t)}{1-re^{i\theta} t}\right)\overline{g(e^{i\theta})}\,d\theta
\\ = \, &\int_{[0,1)}f(t)\left(\int_0^{2\pi}\frac{\overline{g(e^{i\theta})}d\theta}{1-re^{i\theta}t}\right) d\mu(t)
=\, \int_{[0,1)}f(t)\overline{g(rt)}\,d\mu(t),\quad 0\le
r<1.\nonumber
\end{align}\par\medskip
(i)\,$\Rightarrow $\, (ii). Assume that $\nu $ is a Carleson measure
and take $f\in \mathcal B$ and $g\in H^1$. Using (\ref{duality1})
and (\ref{gr-Bloch}), we obtain \begin{align*}& \left|\int_0^{2\pi}
I_\mu(f)(re^{i\theta})\overline{g(e^{i\theta})}\,d\theta\right|
\,=\, \left|\int_{[0,1)}f(t)\overline{g(rt)} \,d\mu(t)\right|\\
& \lesssim \, \|f\|_{\mathcal B}
\int_{[0,1)}\left|g(rt)\right|\log\frac{2}{1-t}\,d\mu(t)\,=\,
\|f\|_{\mathcal B}
\int_{[0,1)}\left|g(rt)\right|\,d\nu(t).\end{align*} Since $\nu $ is
a Carleson measure
$$\int_{[0,1)}\left|g(rt)\right|\,d\nu (t)\,\lesssim
\,\|g_r\|_{H^1}\,\le \,\|g\|_{H^1}.$$ Here, as usual, $g_r$ is the
function defined by $g_r(z)=g(rz)$ ($z\in \mathbb D$).
\par
Thus, we have proved that
$$\left|\int_0^{2\pi}
I_\mu(f)(re^{i\theta})\overline{g(e^{i\theta})}\,d\theta\right|
\,\lesssim \,\|f\|_{\mathcal B}\|g\|_{H^1},\quad f\in \mathcal
B,\,\, g\in H^1.$$ Using Fefferman's duality Theorem (see
\cite[Theorem\,\@7.\,\@1]{G:BMOA}) we deduce that if $f\in \mathcal
B$ then $I_\mu (f)\in BMOA$ and
$$\Vert I_\mu (f)\Vert_{BMOA}\lesssim \Vert f\Vert _{\mathcal B}.$$
\par\bigskip The implication (ii)\, $\Rightarrow $ \, (iii) is
trivial because $BMOA\subset \mathcal B$.
\par\bigskip (iii)\, $\Rightarrow $ \, (i).
Assume (iii). Then there exists a positive constant $A$ such that
$\Vert I_\mu (f)\Vert _{BMOA}\le A\Vert f\Vert _{BMOA}$, for all
$f\in BMOA$. Set
$$F(z)=\log\frac{2}{1-z},\quad z\in \mathbb D.$$ It is well known
that $F\in BMOA$. Then $I_\mu (F)\in BMOA$ and $$\Vert I_\mu
(F)\Vert _{BMOA}\le A\Vert F\Vert _{BMOA}.$$ Then using again
Fefferman's duality theorem we obtain that
$$\left\vert \int_{0}^{2\pi }I_\mu (F)(re^{i\theta })\,\overline
{g(e^{i\theta })}\,d\theta \right \vert \lesssim \Vert g\Vert
_{H^1},\quad g\in H^1.$$ Using (\ref{duality1}) and the definition
of $F$, this implies
\begin{equation}\label{Fgbound}\left\vert
\int_{[0,1)]}\,\overline{g(rt)}\,\log\frac{2}{1-t}\,d\mu
(t)\right\vert \,\lesssim\,\Vert g\Vert _{H^1},\quad g\in
H^1.\end{equation}
\par Take $g\in H^1$. Using Proposition\,\@2 of \cite{Ch-Gi-Pe} we know
that there exists a function $G\in H^1$ with $\Vert G\Vert
_{H^1}\,=\,\Vert g\Vert _{H^1}$ and such that
$$\vert g(s)\vert \le G(s),\quad \text{for all $s\in [0,1)$}.$$
Using these properties and (\ref{Fgbound}) for $G$, we obtain
\begin{align*}
\int_{[0,1)}\left|g(rt)\right|\log\frac{2}{1-t}\,d\mu(t)&\le\int_{[0,1)}G(rt)\log\frac{2}{1-t}\,d\mu(t)\\
&\le C \|G_r\|_{H^1}\le C \|G\|_{H^1}=C\|g\|_{H^1}
\end{align*}
for a certain constant $C>0$, independent of $g$. Letting $r$ tend
to $1$, it follows that
$$ \int_{[0,1)}\left|g(t)\right|\log\frac{2}{1-t}\,d\mu(t)\lesssim
\|g\|_{H^1},\quad g\in H^1.$$ This is equivalent to saying that $\nu
$ is a Carleson measure.
\end{Pf}
\par\bigskip
It is worth noticing that for $\mu $ and $\nu $ as in
Theorem~\ref{ImuBBMOA}, $\nu $ being a Carleson measure is
equivalent to $\mu $ being an $1$-logarithmic $1$-Carleson measure.
Actually, we have the following more general result.
\begin{proposition}\label{mu-nu} Let $\mu $ be a positive Borel
measure on $[0, 1)$, $s>0$, and $\alpha \ge 0$. Let $\nu $ be the
Borel measure on $[0, 1)$ defined by $$d\nu (t)\,=\,\left
(\log\frac{2}{1-t}\right )^\alpha d\mu (t).$$ Then, the following
two conditions are equivalent.
\begin{itemize}\item [(a)] $\nu $ is an $s$-Carleson measure.
\item [(b)] $\mu $ is an $\alpha $-logarithmic $s$-Carleson measure.
\end{itemize}
\end{proposition}
\par \begin{pf}
\par (a)\,$\Rightarrow $\, (b). Assume (a). Then there exists a
positive constant $C$ such that
$$\int_{[t, 1)}\left (\log \frac{2}{1-u}\right )^{\alpha }\,d\mu (u)\,\le \,C(1-t)^s,\quad t\in
[0,1).$$ Using this and the fact that the function $u\,\mapsto
\,\log \frac{2}{1-u}$ is increasing in $[0, 1)$, we obtain
\begin{align*}\left (\log \frac{2}{1-t}\right )^\alpha \int_{[t,
1)}d\mu (u)\,\le \,\int_{[t, 1)}\left (\log \frac{2}{1-u}\right
)^\alpha d\mu (u)\,\le \,C(1-t)^s,\quad t\in [0,1). \end{align*}
This shows that $\mu $ is an $\alpha $-logarithmic $s$-Carleson
measure.
\par (b)\,$\Rightarrow $\, (a). Assume (b). Then there exists a
positive constant $C$ such that \begin{equation}\label{mulog}\left
(\log \frac{2}{1-t}\right )^\alpha  \mu \left ([t, 1)\right
)\,\le\,C\,(1-t)^s,\quad 0\le t<1.\end{equation} For $0\le u<1$, set
$F(u)=\mu \left ([0,u)\right )\,-\,\mu \left ([0,1)\right )\,=\,-\mu
\left ([u,1)\right )$. Integrating by parts and using (\ref{mulog}),
we obtain
\begin{align*}& \nu \left ([t,1)\right )\,=\,
\int_{[t,1)}\left (\log\frac{2}{1-u}\right )^\alpha \,d\mu (u)
\\
=&\,\left (\log \frac{2}{1-t}\right )^\alpha \mu \left ([t,1)\right
)\,-\,\lim_{u\to 1^-}\left (\log\frac{2}{1-u}\right )^\alpha \mu
\left ([u,1)\right )\\&+\,\alpha \int_{[t,1)}\mu \left ([u,1)\right
)\left (\log\frac{2}{1-u}\right )^{\alpha -1}\frac{du}{1-u}
\\
=&\,\left (\log\frac{2}{1-t}\right )^\alpha \mu \left ([t,1)\right
)\,+\,\alpha \int_{[t,1)}\mu \left ([u,1)\right )\left
(\log\frac{2}{1-u}\right )^{\alpha -1}\frac{du}{1-u}\\
\le &\, C\,(1-t)^s \,+\,C\,\alpha
\,\int_t^1\,\frac{(1-u)^{s-1}}{\log \frac{2}{1-u}}\,du\\\lesssim
&\,(1-t)^s,\quad 0\le t<1.
 \end{align*} Thus, $\nu $ is an $s$-Carleson measure.
\end{pf}
\par\bigskip
 The following lemma will be needed in the proof of
Theorem\,\@\ref{Imu=Hmu}. \begin{lemma}\label{lem-moments} Let $\mu$
be a positive Borel measure in $[0,1)$ such that the measure $\nu $
defined by $d\nu(t)=\log\frac{1}{1-t}d\mu(t)$ is a Carleson measure.
Then the sequence of moments $\{ \mu _n\} $ satisfies
$$\mu_n=\og \left(\frac{1}{n\log n}\right),\quad \text{as $n\to \infty $}.$$
\end{lemma}
\par\bigskip
Actually, we shall prove the following more general result.
\begin{lemma}\label{gen-lem-moments} Suppose that $0\le \alpha \le \beta $, $s\ge 1$, and let
$\mu$ be a positive Borel measure on $[0,1)$ which is a $\beta
$-logarithmic $s$-Carleson measure. Then
$$\int_{[0,1)}\,t^k\left (\log\frac{2}{1-t}\right )^\alpha \,d\mu
(t)\,=\,\og \left (\frac{(\log k)^{\alpha -\beta }}{k^s}\right
),\quad\text{as $k\to \infty $}.$$
\end{lemma}
\par\medskip Using Proposition\,\@\ref{mu-nu},
Lemma\,\@\ref{lem-moments} follows taking $\alpha =0$, $\beta =1$,
and $s=1$ in Lemma\,\@\ref{gen-lem-moments}.
\par\medskip

\begin{Pf}{\,\em{Lemma\,\@\ref{gen-lem-moments}.}} Arguing as in the proof of the implication (b)\,$\Rightarrow
$\,(a) of Proposition\,\@\ref{mu-nu}, integrating by parts and using
the fact that $\mu $ is a $\beta $-logarithmic $1$-Carleson measure,
we obtain
\begin{align}\label{est-mom}
& \int_{[0,1)}t^k\left (\log \frac{2}{1-t}\right )^\alpha d\mu (t)
\\ = &\, k\int_{0}^1\mu\big([t,1)\big) t^{k-1} \left (\log
\frac{2}{1-t}\right )^\alpha\,dt \,+\,\alpha
\int_{0}^1\mu\big([t,1)\big) t^{k}\left (\log \frac{2}{1-t}\right
)^{\alpha -1}\frac{dt}{1-t} \nonumber
\\ \lesssim &\,
k\int_{0}^1 (1-t)^st^{k-1} \left (\log \frac{2}{1-t}\right )^{\alpha
-\beta}\,dt \,+\,\alpha \int_{0}^1(1-t)^{s-1}t^{k}\left (\log
\frac{2}{1-t}\right )^{\alpha -\beta -1}{dt}. \nonumber
\end{align}
Now, we notice that the weight functions $$\omega_1 (t)=(1-t)^s\left
(\log\frac{2}{1-t}\right )^{\alpha -\beta }\,\,\,\text{and}\quad
\omega_2 (t)=(1-t)^{s-1}\left (\log\frac{2}{1-t}\right )^{\alpha
-\beta -1}$$ are regular in the sense of \cite{PeRa-Mem} (see
\cite[p.\,\@6]{PeRa-Mem} and \cite[Example\,\@2]{AlSi}). Then, using
Lemma\,\@1.\,\@3 of \cite{PeRa-Mem} and the fact that the
$\omega_j$'s are also decreasing, we obtain \begin{align*}\int_{0}^1
(1-t)^st^{k-1} \left (\log \frac{2}{1-t}\right )^{\alpha -\beta}dt\,
\lesssim &\,\int_{1-\frac{1}{k}}^1 (1-t)^st^{k-1} \left (\log
\frac{2}{1-t}\right )^{\alpha -\beta}dt\\ \lesssim &\,\frac{(\log
k)^{\alpha -\beta}}{k^{s+1}}\end{align*} and
\begin{align*}\int_{0}^1(1-t)^{s-1}t^{k}\left (\log \frac{2}{1-t}\right
)^{\alpha -\beta -1}{dt}\,\lesssim &\,
\int_{1-\frac{1}{n}}^1(1-t)^{s-1}t^k\left (\log \frac{2}{1-t}\right
)^{\alpha -\beta -1}{dt}\\\lesssim &\,\frac{(\log k)^{\alpha -\beta
-1}}{k^s}.\end{align*}
 Using these two estimates in (\ref{est-mom}) yields
$$\int_{[0,1)}t^k\left (\log \frac{2}{1-t}\right )^\alpha d\mu (t)\,\lesssim \, \frac{(\log k)^{\alpha -\beta }}{k^s}$$ finishing the proof.
\end{Pf}
\par\bigskip
We shall also use the characterization of the coefficient
multipliers from $\mathcal B$ into $\ell^1$ obtained by Anderson and
Shields in \cite{AnSh}.
\begin{other}\label{mult-Bell1} A sequence $\{ \lambda _n\}
_{n=0}^\infty $ of complex numbers is a coefficient multiplier from
$\mathcal B$ into $\ell ^1$ if and only if
\begin{equation*}\label{mult-cond}\sum_{n=1}^\infty\left( \sum_{k=2^n+1}^{2^{n+1}}
|\lambda_k|^2\right)^{1/2}<\infty.\end{equation*}
\end{other}
Bearing in mind Definition\,\@1 of \cite{AnSh}, Theorem\,\@A reduces
to the case $p=1$ in Corollary\,\@1 in p.\,\@259 of \cite{AnSh}.
\par We recall that if $X$ is a space of analytic
functions in $\mathbb D$ and $Y$ is a space of complex sequences, a
sequence $\{ \lambda _n\} _{n=0}^\infty \subset \mathbb C$ is said
to be a multiplier of $X$ into $Y$ if whenever
$f(z)=\sum_{n=0}^\infty a_nz^n\in X$ one has that the sequence $\{
\lambda _na_n\} _{n=0}^\infty $ belongs to $Y$. Thus:
\par By saying that $\{ \lambda_n \}_{n=0}^\infty $ is a coefficient multiplier from $\mathcal B$ into $\ell ^1$ we
mean that  $$\text{If\, $f(z)=\sum_{n=0}^\infty a_nz^n\in \mathcal
B$\,\, then \,\, $\sum_{n=0}^\infty \vert \lambda _na_n\vert <\infty
$.}$$ Actually, using the closed graph theorem, we can assert the
following:
\par A complex sequence $\{ \lambda _n\} _{n=0}^\infty $ is a
multiplier from $\mathcal B$ to $\ell ^1$ if and only if there
exists a positive constant $C$ such that whenever
$f(z)=\sum_{n=0}^\infty a_nz^n\in \mathcal B$, we have that
$\sum_{n=0}^\infty \vert \lambda _na_n\vert \le C\Vert f\Vert
_{\mathcal B}.$\par\bigskip
\begin{Pf}{\,\em{Theorem \ref{Imu=Hmu}.}}
Suppose that $\nu $ is a Carleson measure. Then, using
Lemma~\,\@\ref{lem-moments}, we see that there exists $C>0$ such
that
\begin{equation}\label{mun}\vert \mu _n\vert \le \frac{C}{n\log
n},\quad n\ge 2.\end{equation} It is clear that $$k^2\log^2k\ge
2^{2n}n^2(\log 2)^2,\quad \text{if\, $2^n+1\le k\le 2^{n+1}$\, for
all \, $n$.}$$ Then it follows that
\begin{align*}
\sum_{n=1}^\infty\left( \sum_{k=2^n+1}^{2^{n+1}} \frac{1}{k^2\log^2
k}\right)^{1/2}\lesssim \sum_{n=1}^\infty\left(\frac{2^n}{n^2
2^{2n}} \right)^{1/2}=\sum_{n=1}^\infty \frac{1}{n 2^{n/2}}<\infty.
\end{align*}
Using this, (\ref{mun}) and Theorem~\,\@\ref{mult-Bell1}, we obtain:
\begin{equation}\label{seq-mom}\text{The sequence of moments $\{ \mu_n\}_{n=0}^\infty $ is a multiplier
from $\mathcal B$ to $\ell^1$.}\end{equation}
\par Take now $f\in \mathcal B$, $f(z)=\sum_{n=0}^\infty a_nz^n$
($z\in \mathbb D$). Using the simple fact that the sequence $\{
\mu_n\}_{n=0}^\infty $ is a decreasing sequence of positive numbers
and (\ref{seq-mom}), we see that there exists $C>0$ such that
\begin{equation}\label{abs-sum}\sum_{k=0}^\infty \vert \mu_{n+k}\,a_k\vert \le \sum_{k=0}^\infty
\vert \mu_k\,a_k\vert \le C\Vert f\Vert _{\mathcal B}, \quad
n\,=\,0,\,1,\,2,\dots .\end{equation} This implies that $\mathcal
H_{\mu }(f)(z)$ is well defined for all $z\in \mathbb D$ and that,
in fact, $\mathcal H_{\mu }(f)$ is an analytic function in $\mathbb
D$. Furthermore, since (\ref{abs-sum}) also implies that we can
interchange the order of summation in the expression defining
$\mathcal H_{\mu }(f)(z)$, we have \begin{align*} \mathcal H_{\mu
}(f)(z)\,= & \,\sum_{n=0}^{\infty}\left(\sum_{k=0}^{\infty}
\mu_{n+k}{a_k}\right)z^n\,=\,\sum_{k=0}^{\infty}a_k\,\left
(\sum_{n=0}^{\infty}\mu_{n+k}z^n\right ) \\ = &\,
\sum_{k=0}^{\infty}a_k\,\left
(\sum_{n=0}^{\infty}\int_{[0,1)}t^{n+k}z^n\,d\mu (t)\right
)=\,\sum_{k=0}^{\infty}\int_{[0,1)}\frac{a_kt^k}{1-tz}\,d\mu (t)\\
=\, & \int_{[0,1)}\frac{f(t)}{1-tz}\,d\mu (t)\,=\,I_\mu (f)(z),\quad
z\in \mathbb D.
\end{align*}
\end{Pf}
\par\bigskip

\par\bigskip
We have the following result regarding compactness.
\begin{theorem}\label{Imucompact} Let $\mu $ be a positive Borel measure on $[0,1)$
with $\int_{[0,1)}\log \frac{2}{1-t}\,d\mu (t)\,<\,\infty $. If the
measure $\nu $ defined by $d\nu (t)=\log \frac{2}{1-t}\,d\mu (t)$ is
a vanishing Carleson measure then:
\begin{itemize}
\item[(i)] The operator $I_\mu $ is a compact operator from $\mathcal B$ into
$BMOA$.
\item[(ii)] The operator $I_\mu $ is a compact operator from $BMOA$ into
itself.
\end{itemize}
\end{theorem}
\par\bigskip
Before embarking on the proof of Theorem\,\@\ref{Imucompact} it is
convenient to recall some facts about Carleson measures and to fix
some notation.
\par If $\mu $ is a Carleson measure on $\mathbb D$, we define the
Carleson-norm of $\mu $, denoted $\mathcal N(\mu )$, as
$$\mathcal N(\mu )\,=\,\sup_{I\, \text{subarc of } \partial \mathbb D}\frac{\mu
\left (S(I)\right )}{\vert I\vert }.$$ We let also $\mathcal E(\mu
)$ denote the norm of the inclusion operator $i:H^1\rightarrow
L^1(d\mu )$. It turns out that these quantities are equivalent:
There exist two positive constants $A_1, A_2$ such that
$$A_1\,\mathcal N(\mu )\,\le \,\mathcal E(\mu )\,\le A_2\mathcal N(\mu
),\quad \text{for every Carleson measure $\mu $ on $\mathbb D$.}$$
\par For a Carleson measure $\mu $ on $\mathbb D$ and $0<r<1$, we
let $\mu _r$ be the measure on $\mathbb D$ defined by
$$d\mu _r(z)\,=\,\chi_{\{ r<\vert z\vert <1\} }d\mu (z).$$
We have that $\mu $ is a vanishing Carleson measure if and only if
$$\mathcal N(\mu _r)\to 0,\quad \text{as $r\to 1$}.$$
\par\bigskip
\par\bigskip
\begin{Pf}{\,\em{Theorem \ref{Imucompact}.}} Since $BMOA$ is
continuously contained in the Bloch spaces, it suffices to prove
(i).
\par Suppose that $\nu $ is a vanishing Carleson
measure. Let $\{ f_n\}_{n=1}^\infty $ be a sequence of Bloch
functions with $\sup_{n\ge 1}\Vert f_n\Vert _{\mathcal B}<\infty $
and such that $\{ f_n\}\to 0$, uniformly on compact subsets of
$\mathbb D$. We have to prove that $I_{\mu }(f_n)\,\to \,0$ in
$BMOA$.
\par The condition $\sup_{n\ge 1}\Vert f_n\Vert
_{\mathcal B}<\infty $ implies that there exists a positive constant
$M$ such that
\begin{equation}\label{boundfn}\vert f_n(z)\vert \,\le
M\log\frac{2}{1-\vert z\vert },\quad z\in \mathbb D,\quad n\ge
1.\end{equation} Recall that for $0<r<1$, $\nu_r$ is the measure
defined by
$$d\nu _r(t)\,=\,\chi _{\{ r<t<1\} }\,d\nu (t).$$
Since $\nu $ is a vanishing Carleson measure, we have that $\mathcal
N(\nu _r)\to 0$, as $r\to 1$, or, equivalently,
\begin{equation}\label{Erto0}\mathcal E(\nu _r)\,\to \,0,\quad
\text{as $t\,\to \,1$}.\end{equation} Take $g\in H^1$ and $r\in
[0,1)$. Using (\ref{boundfn}) we have
\begin{align*}\int_{[0,1)}\vert f_n(t)\vert \vert g(t)\vert \,d\mu
(t)\,=\,&\int_{[0,r)}\vert f_n(t)\vert \vert g(t)\vert \,d\mu
(t)\,+\,\int_{[r,1)}\vert f_n(t)\vert \vert g(t)\vert \,d\mu (t)
\\
\le \,&\int_{[0,r)}\vert f_n(t)\vert \vert g(t)\vert \,d\mu
(t)\,+\,M\int_{[r,1)}\log\frac{2}{1-t}\vert g(t)\vert \,d\mu (t)\\
=\,& \int_{[0,r)}\vert f_n(t)\vert \vert g(t)\vert \,d\mu
(t)\,+\,M\int_{[0,1)}\vert g(t)\vert \,d\nu_r (t)\\ \le \,&
\int_{[0,r)}\vert f_n(t)\vert \vert g(t)\vert \,d\mu
(t)\,+\,M\mathcal E(\nu _r)\Vert g\Vert _{H^1}.
\end{align*}
Using (\ref{Erto0}) and the fact that $\{ f_n\}\to 0$, uniformly on
compact subsets of $\mathbb D$, it follows that
$$\lim_{n\to \infty }\int_{[0,1)}\vert f_n(t)\vert \vert g(t)\vert \,d\mu
(t)\,=\,0,\quad\text{for all $g\in H^1$}.$$ Bearing in mind
(\ref{duality1}), this yields
$$\lim_{n\to\infty }\left (\lim_{r\to 1}\left \vert \int_0^{2\pi
}\,I_\mu (f_n)(re^{i\theta })\overline {g(e^{i\theta })}\,d\theta
\right \vert \right )\,=\, 0,\quad\text{for all $g\in H^1$}.$$ By
the duality relation $\left (H^1\right)^\star \,=\,BMOA$, this is
equivalent to saying that $I_{\mu }(f_n)\,\to \,0$ in $BMOA$.
\end{Pf}
\par\bigskip
\section{The operator $\mathcal H_\mu $ acting on $Q_s$ spaces and Besov spaces}\label{section-Besov-Qs}
\par If $0\le s<\infty$, we say that $f\in Q_s$ if $f$ is analytic in
$\D$ and \begin{equation*}\Vert f\Vert_{Q_s}\,\defeq \,\left (\vert
f(0)\vert ^2\,+\,\rho_{Q_s}(f)^2\right )^{1/2}\,<\,\infty
,\end{equation*} where
\[
\rho_{Q_s}(f)\,\defeq \,\left (
\sup_{a\in\D}\int_\D|f'(z)|^2g(z,a)^s\,dA(z)\right )^{1/2}.\] Here,
$g(z,a)$ is the Green's function in $\D$, given by
$g(z,a)=\log\left|\frac{1-\overline{a}z}{z-a}\right|$, while
$dA(z)=\frac{dx\,dy}{\pi }$ is the normalized area measure on
$\mathbb D$. All $Q_s$ spaces ($0\le s<\infty $) are conformally
invariant with respect to the semi-norm $\rho _{Q_s}$ (see e.\@g.,
\cite[p.\,@1]{X2} or \cite[p.\,\@47]{DGV1}).
\par
These spaces were introduced by Aulaskari and Lappan
in~\cite{au-la94} while looking for new characterizations of Bloch
functions. They proved that for $s>1$, $Q_s$ is the Bloch space.
Using one of the many characterizations of the space $BMOA$ (see,
e.\,\@g.,  \cite[Theorem\,\@5]{Ba} or
\cite[Theorem\,\@6.\,\@2]{G:BMOA}) we see that $Q_1=BMOA$. In the
limit case $s=0$, $Q_s$ is the classical Dirichlet space $\mathcal
D$ of those analytic functions $f$ in $\D$ satisfying $\,
\int_\D|f'(z)|^2\,dA(z) <\infty$.
\par
It is well known that $\mathcal D\subset VMOA$. Aulaskari, Xiao and
Zhao proved in~\cite{au-xi-zh95} that
\[\mathcal D\subsetneq Q_{s_1}\subsetneq Q_{s_2}\subsetneq BMOA,\qquad0<s_1<s_2<1.\]
We mention the book \cite{X2} as an excelent reference for the
theory of $Q_s$-spaces.
\par It is well known that the function $F(z)\,=\,\log
\frac{2}{1-z}$ belong to $Q_s$, for all $s>0$, (in fact, it is
proved in \cite{au-la-xi-zh97} that the univalent functions in all
$Q_s$-spaces ($0<s<\infty $) are the same). Using this we easily see
that Theorem\,\@\ref{ImuBBMOA} and Theorem\,\@\ref{HmuboundedBMOA}
can be improved as follows.\par\bigskip
\begin{theorem}\label{ImuBBMOAQs} Let $\mu $
be a positive Borel measure on $[0,1)$. Then the following
conditions are equivalent:
\begin{itemize}
\item[(i)] $\int_{[0,1)}\log \frac{2}{1-t}\,d\mu (t)\,<\,\infty $.
\item[(ii)] For any given $s\in (0, \infty )$ and any $f\in Q_s$, the integral in (\ref{Imu})
converges for all $z\in \mathbb D$ and the resulting function $I_\mu
(f)$ is analytic in \,$\D $.
\end{itemize}
We remark that condition (ii) with $s\ge 1$ includes the points (ii)
and (iii) of Theorem\,\@\ref{ImuBBMOA}.
\end{theorem}
\par\bigskip
\begin{theorem}\label{ImuboundedBMOAQs} Let $\mu $ be a positive Borel measure on $[0,1)$
with $\int_{[0,1)}\log \frac{2}{1-t}\,d\mu (t)\,<\,\infty $. Then
the following two conditions are equivalent:
\begin{itemize}\item[(i)] The measure $\nu $ defined by $d\nu
(t)=\log \frac{2}{1-t}\,d\mu (t)$ is a Carleson measure.
\item[(ii)] For any given $s\in (0, \infty )$, the  operator $I_\mu $ is bounded from $Q_s$ into
$BMOA$.
\end{itemize}
We remark that (ii) with $s>1$ reduces to condition (ii) of
Theorem\,\@\ref{ImuboundedBMOA}, while (ii) with $s=1$ reduces to
condition (iii) of Theorem\,\@\ref{ImuboundedBMOA}.
\end{theorem}
\par\bigskip
These results cannot be extended to the limit case $s=0$. Indeed,
the function $F(z)\,=\,\log\frac{2}{1-z}$ does not belong to the
Dirichlet space $\mathcal D$.
\par\bigskip The Dirichlet space is one among the analytic Besov spaces.
\par For $1<p<\infty $, the {\it analytic Besov space\/} $B\sp p $ is
defined as the set of all functions $f $ analytic in $\mathbb D $
such that $$\Vert f\Vert_{B^p}\defeq \left (\vert f(0)\vert
^p\,+\,\rho_ p(f)^p\right )^{1/p}<\infty ,$$ where
$$\rho_ p(f)=\left (\int_{\mathbb D}(1-\vert z\vert^2)^{p-2}\vert
f^\prime (z)\vert ^p\,dA(z)\right )^{1/p}.$$ All $B\sp p $ spaces
($1<p<\infty $) are conformally invariant with respect to the
semi-norm $\rho _p$ (see \cite[p.\,\@112]{AFP} or
\cite[p.\,\@46]{DGV1}). We have that $\mathcal D\,=\,B^2$. A lot of
information on Besov spaces can be found in \cite{AFP, DGV1, HW1, Z,
Zhu-book}. Let us recall that
$$B^p\,\subsetneq \,B^q\,\subsetneq VMOA,\quad
1\,<\,p\,<\,q\,<\infty  .$$
\par From now on, if $1\,<\,p\,<\infty $ we let $p^\prime $ denote the exponent conjugate to $p$, that
is, $p^\prime $ is defined by the relation
$\frac{1}{p}\,+\,\frac{1}{p^\prime }\,=\,1$. If $f\in B^p$
($1\,<\,p\,<\,\infty $) then, see \cite{HW1} or \cite{Z},
\begin{equation}\label{grBpop}
\vert f(z)\vert \,=\,\op \left (\left (\log \frac{1}{1-\vert z\vert
}\right )^{1/p^\prime }\right ),\quad \text{as $\vert z\vert \to
1$},\end{equation} and there exists a positive constant $C\,>,0$
such that
\begin{equation}\label{grBp}
\vert f(z)\vert \,\le\,C\Vert f\Vert_{B^p}\left (\log
\frac{2}{1-\vert z\vert }\right )^{1/p^\prime },\quad z\in \mathbb
D,\quad f\in B^p.\end{equation}
\par\medskip
Clearly, (\ref{grBpop}) or (\ref{grBp}) imply that the function
$F(z)\,=\,\log\frac{2}{1-z}$ does not belong to $B^p$ ($1\,<\,p\,<\,
\infty $), a fact that we have already mentioned for $p=2$. Our
substitutes of Theorem\,\@\ref{ImuBBMOA} and
Theorem\,\@\ref{ImuboundedBMOA} for Besov spaces are the following.
\par\medskip
\begin{theorem}\label{ImuBp} Let $1<p<\infty $ and
let $\mu $ be a positive Borel measure on $[0,1)$. We have:
\begin{itemize}
\item[(i)] If $\int_{[0,1)}\left (\log\frac{2}{1-t}\right )^{1/p^\prime }\,d\mu (t)\,<\,\infty
$, then for any given $f\in B^p$, the integral in (\ref{Imu})
converges for all $z\in \mathbb D$ and the resulting function $I_\mu
(f)$ is analytic in \,$\D $.
\item[(ii)] If for any given $f\in B^p$, the integral in (\ref{Imu})
converges for all $z\in \mathbb D$ and the resulting function $I_\mu
(f)$ is analytic in \,$\D $, then $\int_{[0,1)}\left(\log
\frac{2}{1-t}\right )^{\gamma }\,d\mu (t)\,<\,\infty $ for all
$\gamma \, <\,\frac{1}{p^\prime }.$
\end{itemize}
\end{theorem}
\par\medskip
\begin{theorem}\label{ImuboundedBpBMOA} Suppose that $1<p<\infty $ and
let $\mu $ be a positive Borel measure on $[0,1)$. Let $\nu $ be the
measure defined by
$$d\nu (t)\,=\,\left (\log\frac{2}{1-t}\right )^{1/p^\prime }\,d\mu
(t).$$ \begin{itemize}\item[(i)] If $\nu $ is a Carleson measure,
then the operator $I_\mu $ is bounded from $B^p$ into $BMOA$.
\item[(ii)] If $\nu $ is a vanishing Carleson measure then the
operator $I_\mu $ is compact from $B^p$ into $BMOA$.
\end{itemize}
\end{theorem}
\par\medskip
These results follow using  the growth condition (\ref{grBp}), the
fact that if $\gamma \,<\,\frac{1}{p^\prime }$ then the function
$f(z)=\left (\log \frac{2}{1-z}\right )^\gamma $ belongs to $B^p$
(see \cite[Theorem\,\@1]{HW1}), and with arguments similar to those
used in the proofs of Theorem\,\@\ref{ImuBBMOA},
Theorem\,\@\ref{ImuboundedBMOA}, and Theorem\,\@\ref{Imucompact}. We
omit the details.
\par\bigskip
Let us work next with the operator $\mathcal H_\mu $ directly. In
order to study its action on the Besov spaces we need some results
on the Taylor coefficients of functions in $B^p$. The following
result was proved by Holland and Walsh in \cite[Theorem~2]{HW1}.
\begin{other}\label{HW-coef-Bp} \par\begin{itemize}
\item[(i)] Suppose that $1<p\le 2$. Then there exists a positive
constant $C_p$ such that if $f\in B^p$ and $f(z)=\sum_{k=0}^\infty
a_k\,z^k$ $(z\in \mathbb D)$ then
\begin{equation*}\sum_{k=1}^\infty k^{p-1}\vert a_k\vert
^p\,\le \,C_p\,\rho_p(f)^p.\end{equation*}
\item[(ii)] If $2\,\le p\,<\infty $ then there exists $C_p>0$ such that if $f(z)=\sum_{k=0}^\infty
a_k\,z^k$ $(z\in \mathbb D)$ with $\sum_{k=1}^\infty k^{p-1}\vert
a_k\vert ^p\,<\,\infty $ then $f\in B^p$ and
\begin{equation*}\rho_p(f)^p\,\le \,C_p\,\sum_{k=1}^\infty k^{p-1}\vert a_k\vert
^p.\end{equation*}
\end{itemize}
\par If $p\neq 2$ the converses to (i) and (ii) are false.
\end{other}
\par Theorem\,\@\ref{HW-coef-Bp} is the analogue for Besov spaces of
results of Hardy and Littlewood for Hardy spaces (Theorem\,\@6.\,\@2
and Theorem\,\@6.\,\@3 of \cite{D}).
\par In spite of the fact that the converse to (ii) is not true,
the membership of $f$ in $B^p$ ($p>2$) implies some summability
conditions on the Taylor coefficients $\{ a_k\} $ of $f$. Indeed,
Pavlovi\'c has proved the following result in
\cite[Theorem\,\@2.\,\@3]{Pav-Inv-Besov}.
\begin{other}\label{pma2Pavlovic} Suppose that $2\,<\,p\,<\,\infty
$. Then there exists a positive constant $C_p$ such that if $f\in
B^p$ and $f(z)=\sum_{k=0}^\infty a_k\,z^k$ $(z\in \mathbb D)$ then
\begin{equation*}\sum_{k=1}^\infty k\vert a_k\vert ^p\,\le C_p\,
\rho_p(f)^p.\end{equation*}
\end{other}
\par\medskip
These results allow us to obtain conditions on $\mu $ which are
sufficient to ensure that $\mathcal H_\mu $ is well defined on the
Besov spaces.
\begin{theorem}\label{Hmu-def-Bp} Let $\mu $ be a finite positive Borel
measure on $[0,1)$.
\begin{itemize}\item[(i)] If $1\,<p\,\le \, 2$ and $\sum_{k=1}^\infty
\frac{\mu_k^{p^\prime }}{k}\,<\,\infty $, then the operator
$\mathcal H_\mu $ is well defined in $B^p$. \item[(ii)] If
$2\,<p\,<\,\infty $ and $\sum_{k=1}^\infty \frac{\mu_k^{p^\prime
}}{k^{p^\prime /p}}\,<\,\infty $, then the operator $\mathcal H_\mu
$ is well defined in $B^p$.
\end{itemize}
\end{theorem}
\begin{pf} Suppose that $1<p<\infty $ and $f\in B^p$,
$f(z)\,=\,\sum_{k=0}^\infty a_kz^k$ ($z\in \mathbb D$). Since the
sequence of moments $\{ \mu _n\} _{n=0}^\infty $ is clearly
decreasing we have
\begin{equation*}\label{coefH-decr}\sum_{k=1}^\infty \vert \mu
_{n+k}\vert \vert a_k\vert \,\le \,\sum_{k=1}^\infty \vert \mu
_{k}\vert \vert a_k\vert ,\quad\text{for all $n\ge
0$}.\end{equation*} Consequently, we have:
\begin{itemize}\item[(i)]
\par If $1\,<\,p\le 2$ and $f\in B^p$, $f(z)\,=\,\sum_{k=0}^\infty a_kz^k$ ($z\in \mathbb
D$), then \begin{align*}\sum_{k=1}^\infty \vert
\mu_{n+k}a_k\vert\,\le\,\sum_{k=1}^\infty \vert \mu _{k}\vert \vert
a_k\vert \,=\,\sum_{k=1}^\infty k^{1-\frac{1}{p}}\vert a_k\vert
\frac{\mu_k}{k^{1/p^\prime }},\quad n\ge 0.\end{align*} Then using
H\"{o}lder inequality and Theorem\,\@\ref{HW-coef-Bp}\,\@(i), we
obtain
\begin{align*}\sum_{k=1}^\infty \vert \mu_{n+k}a_k\vert\,\le &\,\left
(\sum_{k=1}^\infty k^{p-1}\vert a_k\vert^p\right )^{1/p}\left
(\sum_{k=1}^\infty \frac{\vert \mu _k\vert^{p^\prime }}{k}\right
)^{1/p^\prime }\\ \le &\,C\,\rho_p(f)\left (\sum_{k=1}^\infty
\frac{\vert \mu _k\vert^{p^\prime }}{k}\right )^{1/p^\prime },\quad
n\ge 0.\end{align*} Then it is clear that the condition
$\sum_{k=1}^\infty \frac{\vert \mu _k\vert^{p^\prime
}}{k}\,<\,\infty $ implies that the power series appearing in the
definition of $\mathcal H_\mu (f)$ defines an analytic function in
$\mathbb D$. \item[(ii)]
\par If $2\,<\,p\,<\,\infty $ and $f\in B^p$, $f(z)\,=\,\sum_{k=0}^\infty a_kz^k$ ($z\in \mathbb
D$), then \begin{align*}\sum_{k=1}^\infty \vert
\mu_{n+k}a_k\vert\,\le \,\sum_{k=1}^\infty \vert \mu _{k}\vert \vert
a_k\vert \,=\,\sum_{k=1}^\infty k^{\frac{1}{p}}\vert a_k\vert
\frac{\mu_k}{k^{1/p}},\quad n\ge 0.\end{align*} Then using
H\"{o}lder inequality and Theorem\,\@\ref{HW-coef-Bp}\,\@(ii), we
obtain
\begin{align*}\sum_{k=1}^\infty \vert \mu_{n+k}a_k\vert\,\le &\,\left
(\sum_{k=1}^\infty k\vert a_k\vert^p\right )^{1/p}\left
(\sum_{k=1}^\infty \frac{\vert \mu _k\vert^{p^\prime }}{k^{p^\prime
/p}}\right )^{1/p^\prime } \\\le &\,C\,\rho_p(f)\left
(\sum_{k=1}^\infty \frac{\vert \mu _k\vert^{p^\prime }}{k^{p^\prime
/p}}\right )^{1/p^\prime },\quad n\ge 0.
\end{align*} Then we see that the
condition $\sum_{k=1}^\infty \frac{\vert \mu _k\vert^{p^\prime
}}{k^{p^\prime /p}}\,<\,\infty $ implies that the power series
appearing in the definition of $\mathcal H_\mu (f)$ defines an
analytic function in $\mathbb D$.
\end{itemize}
\end{pf}
\par\bigskip Let us turn to study when is
the operator $\mathcal H_\mu $  bounded from $B^p$ into itself. Let
us mention that Bao and Wulan \cite{Bao-Wu} considered an operator
which is closely related to the operator $\mathcal H_\mu $ acting on
the Dirichlet spaces $\mathcal D_\alpha $ ($\alpha \in \mathbb R$)
which are defined as follows:
\par For $\alpha \in \mathbb R$, the space $\mathcal
D_\alpha $ consists of those functions $f(z)=\sum_{n=0}^\infty
a_n\,z^n$ analytic in $\mathbb D$ for which
$$\Vert f\Vert _{\mathcal D_\alpha }\,\defeq \,\left (\sum_{n=0}^\infty
(n+1)^{1-\alpha }\,\vert a_n\vert ^2\right )^{1/2}\,<\,\infty .$$
Let us remark that $\mathcal D_0$ is the Dirichlet spaces $\mathcal
D\,=\,B^2$, while $\mathcal D_1\,=\,H^2$.
\par Bao and Wulan proved that if $\mu $ is a positive Borel measure on $[0,1)$ and $0<\alpha <2$, then the operator
$\mathcal H_\mu $ is bounded from $\mathcal D_\alpha $ into itself
if and only if $\mu $ is a Carleson measure. Let us remark that this
does not include the case $\alpha =0$. In fact, the following
results are proved in \cite{Bao-Wu}.
\begin{other}\label{Bao-Wulan-Dirichlet} $ $ \par
\begin{itemize}
\item[(i)] There exists a positive Borel measure $\mu $ on $[0,1)$
which is a Carleson measure but such that $\mathcal H_\mu
(B^2)\not\subset B^2$.
\item[(ii)]
Let $\mu $ be a positive Borel measure on $[0,1)$ such that the
operator $\mathcal H_\mu $ is a bounded operator from $B^2$ into
itself. Then $\mu $ is a Carleson measure.
\end{itemize}
\end{other}
\par\medskip We can improve these results and, even more, we shall
obtain extensions of these improvements to all $B^p$ spaces
($1<p<\infty $). More precisely we are going to prove the following
results.

\begin{theorem}\label{Bp-not-suf} Suppose that $1<p<\infty $ and $0<\beta \le \frac{1}{p}$. Then
there exists a positive Borel measure $\mu $ on $[0,1)$ which is a
$\beta $-logarithmic $1$-Carleson measure but such that the operator
$\mathcal H_\mu $ does not apply $B^p$ into itself.
\end{theorem}
\par\medskip Next we prove that $\mu $ being a $\beta $-logarithmic
$1$-Carleson measure for a certain $\beta $ is a necessary condition
for $\mathcal H_\mu $ being a bounded operator from $B^p$ into
itself.
\begin{theorem}\label{Bp-log-nec} Suppose that $1<p<\infty $ and let
$\mu $ be a positive Borel measure on $[0,1)$ such that the operator
$\mathcal H_\mu $ is bounded from $B^p$ to itself. Then $\mu $ is a
$\gamma $-logarithmic $1$-Carleson measure for any $\gamma
<\,1-\frac{1}{p}$.
\end{theorem}
\par\medskip
Finally, we obtain a sufficient condition for the boundedness of
$\mathcal H_\mu $ from $B^p$ into itself.

\begin{theorem}\label{log-car-suf} Suppose that $1<p<\infty $, $\gamma >1$, and
let $\mu $ be a positive Borel measure on $[0,1)$ which is a $\gamma
$-logarithmic $1$-Carleson measure. Then the operator $\mathcal
H_\mu $ is a bounded operator from $B^p$ into itself.
\end{theorem}
We shall need a number of results on Besov spaces,  as well as some
lemmas, to prove these three theorems. First of all we notice that
the Besov spaces can be characterized in terms of \lq\lq dyadic
blocks\rq\rq . In order to state this in a precise way we need to
introduce some notation. \par For a function $f(z)=\sum_{n=0}^\infty
a_nz^n$ analytic in $\D,$ define the polynomials $\Delta_jf$ as
follows:
\[
\Delta_jf(z)=\sum_{k=2^j}^{2^{j+1}-1} a_kz^k, \quad\text{for $j\ge
1$},
\]
$$\Delta_0f(z)=a_0+a_1z.$$
Mateljevi\'c and Pavlovi\'c proved in
\cite[Theorem\,\@2.\,\@1]{MP-stu} (see also
\cite[Theorem\,\@C]{Pav-dec}) the following result.
\begin{other}\label{block-Bergman}
Let $1<p<\infty $ and $\alpha >-1$.  For a function $f$ analytic in
$\mathbb D$ we define
\begin{align*}Q_1(f)&\,\defeq \, \int_{\mathbb D}\vert f(z)\vert ^p(1-\vert
z\vert )^{\alpha }\,dA(z),\quad Q_2(f)\,\defeq \, \sum_{n=0}^\infty
2^{-n(\alpha +1)}\Vert \Delta _nf\Vert _{H^p}^p.\end{align*} Then,
$Q_1(f)\,\asymp\, Q_2(f)$.
\end{other}
\par\medskip Theorem\,\@\ref{block-Bergman} readily implies the
following result.
\begin{corollary}\label{block-Besov} Suppose that $1<p<\infty $ and
$f$ is an analytic function in $\mathbb D$. Then
$$f\in B^p \,\,\Leftrightarrow \,\,\sum_{n=0}^\infty
2^{-n(p-1)}\Vert \Delta_nf^\prime \Vert _{H^p}^p\,<\,\infty .$$
Furthermore,
$$\rho_p(f)^p\,\asymp \,\sum_{n=0}^\infty
2^{-n(p-1)}\Vert \Delta_nf^\prime \Vert _{H^p}^p.$$
\end{corollary}
\par Using Corollary\,\@\ref{block-Besov} we can prove that the
converses of (i) and (ii) in Theorem\,\@\ref{HW-coef-Bp} hold if the
sequence of Taylor coefficients $\{ a_n\} $ decreases to $0$. This
is the analogue for Besov spaces of the result proved in
\cite{HL-31} by Hardy and Littlewood for Hardy spaces (see also
\cite{Pav-dec}, \cite[7.\,\@5.\,\@9]{Pabook} and
\cite[Chapter\,\@XII, Lemma\,\@6.\,\@6]{Zy}).
\begin{theorem}\label{coef-dec-Besov}
Suppose that $1<p<\infty $ and let $\{ a_n\} _{n=0}^\infty $ be a
decreasing sequence of non-negative numbers with $\{ a_n\} \to 0$,
as $n\to \infty $. Let $f(z)=\sum_{n=0}^\infty a_nz^n$ ($z\in
\mathbb D$). Then
$$f\in B^p\,\,\Leftrightarrow \,\, \sum_{n=1}^\infty n^{p-1}
a_n^p\,<\,\infty .$$ Furthermore, $\rho_p(f)^p\,\asymp
\,\sum_{n=1}^\infty n^{p-1}a_n^p.$
\end{theorem}
\begin{pf} For every $n$,
we have
$$z\left (\Delta_nf^\prime \right
)(z)\,=\,\sum_{k=2^n+1}^{2^{n+1}}ka_kz^k.$$ Since the sequence
$\lambda =\{ k\} _{k=0}^\infty $ is an increasing sequence of
non-negative numbers, using Lemma\,\@A of \cite{Pav-dec} we see that
\begin{equation}\label{Deltanprime}\Vert z\left (\Delta_nf^\prime \right
)\Vert_{H^p}^p\,\asymp \,2^{np}\Vert \Delta_nf\Vert
_{H^p}^p.\end{equation} Now, set $h(z)=\sum_{n=0}^\infty z^n$ ($z\in
\mathbb D$). Since the sequence $\tilde\lambda = \{ a_n\}
_{n=0}^\infty $ is a decreasing sequence of non-negative numbers,
using the second part of Lemma\,\@A of \cite{Pav-dec}, we see that
\begin{equation}\label{Deltanh}a_{2^n}^p\Vert
\Delta_nh\Vert_{H^p}^p\,\lesssim \,\Vert
\Delta_nf\Vert_{H^p}^p\,\lesssim \,a_{2^{n-1}}^p\Vert
\Delta_nh\Vert_{H^p}^p.\end{equation} Notice that
$h(z)=\frac{1}{1-z}$ ($z\in \mathbb D$). Then it is well known that
$M_p(r,h)\,\asymp\,(1-r)^{\frac{1}{p}-1}$ (recall that $1<p<\infty
$). Following the notation of \cite{MP-stu}, this can be written as
$h\in H\left (p,\infty ,1-\frac{1}{p}\right )$. Then using
Theorem\,\@2.\,\@1 of \cite{MP-stu} (see also
\cite[p.\,\@120]{Pabook}), we deduce that $\Vert \Delta_n\Vert_
{H^p}^p\,\asymp\,2^{n(p-1)}$. Using this and (\ref{Deltanh}), it
follows that
\begin{equation}\label{Deltanf}2^{n(p-1)}a_{2^n}^p\,\lesssim\,\Vert
\Delta_nf\Vert_{H^p}^p\,\lesssim
\,2^{n(p-1)}a_{2^{n-1}}^p.\end{equation} Using
Corollary\,\@\ref{block-Besov}, (\ref{Deltanprime}), and
(\ref{Deltanf}), we see that
$$\rho_p(f)^p\,\asymp\,
\sum_{n=0}^\infty 2^{-n(p-1)}\Vert z\,\Delta_nf^\prime \Vert
_{H^p}^p\,\asymp \,\sum_{n=0}^\infty 2^n\Vert
\Delta_nf\Vert_{H^p}^p\,\asymp \,\sum_{n=0}^\infty
2^{np}a_{2^n}^p.$$ Now, the fact that $\{ a_n\} $ is decreasing
implies that $\sum_{n=0}^\infty 2^{np}a_{2^n}^p\,\asymp
\,\sum_{n=1}^\infty n^{p-1}a_n^p$ and, then it follows that
$\rho_p(f)^p\,\asymp\,\sum_{n=1}^\infty n^{p-1}a_n^p.$
\end{pf}
\par\medskip
\begin{remark} If $f$ is an analytic function in $\mathbb D$,
$f(z)=\sum_{n=0}^\infty a_nz^n$ ($z\in \mathbb D$), and $1<p<\infty
$ then any of the two conditions $f\in B^p$ and $\sum_{n=1}^\infty
n^{p-1}\vert a_n\vert ^p\,<\,\infty $ implies that $\{ a_n\} \to 0$.
Consequently, the condition $\{ a_n\} \to 0$ can be omitted in the
hypotheses of Theorem\,\@\ref{coef-dec-Besov}.
\end{remark}
\par\medskip Suppose that $\beta \ge 0$, $s\ge 1$, $1<p<\infty $,
and  $\mu $ is a positive Borel measure on $[0,1)$ which is a $\beta
$-logarithmic $s$-Carleson measure. Using
Lemma\,\@\ref{gen-lem-moments} and Theorem\,\@\ref{Hmu-def-Bp}, it
follows that $\mathcal H_\mu $ is well defined on $B^p$. Also, it is
easy to see that $\int_{[0,1)}\left (\log \frac{2}{1-t}\right
)^{1/p^\prime }\,d\mu (t)\,<\,\infty $, a fact that, using
Theorem\,\@\ref{ImuBp}\,\@(i), shows that $I_\mu $ is also well
defined in $B^p$. Using then standard arguments it follows that
$I_\mu $ and $\mathcal H_\mu $ coincide in $B^p$. Let us state this
as a lemma.
\begin{lemma}\label{HmuImuBp} Suppose that
$\beta \ge 0$, $s\ge 1$, $1<p<\infty $, and  $\mu $ is a positive
Borel measure on $[0,1)$ which is a $\beta $-logarithmic
$s$-Carleson measure. Then the operators $\mathcal H_\mu $ and
$I_\mu $ are well defined in $B^p$ and $\mathcal H_\mu (f)\,=\,I_\mu
(f)$, for all $f\in B^p$.
\end{lemma}

\par\bigskip
\begin{Pf}{\,\em{Theorem \ref{Bp-not-suf}.}}
Let $\mu $ be the Borel measure on $[0,1)$ defined by $$d\mu
(t)\,=\,\left (\log \frac{2}{1-t}\right )^{-\beta }dt.$$ Since the
function $x\,\mapsto \,\left (\log \frac{2}{1-x}\right )^{-\beta }$
is decreasing in $[0,1)$, we have
$$\mu \left ([t,1)\right )\,=\,\int_t^1\left (\log \frac{2}{1-x}\right )^{-\beta
}\,dx\,\le \,(1-t)\left (\log \frac{2}{1-t}\right )^{-\beta },\quad
0\le t<1.$$ Hence, $\mu $ is a $\beta $-logarithmic $1$-Carleson
measure. Then, taking $\alpha =0$ in Lemma\,\@\ref{gen-lem-moments},
we see that $$\mu_k\,=\,\og \left (\frac{1}{k(\log k)^\beta }\right
).$$ On the other hand,
$$\mu_k\,\ge \,\int_0^{1-\frac{1}{k}}t^k\left
(\log\frac{2}{1-t}\right )^{-\beta }\,dt\,\gtrsim \frac{1}{(\log
k)^\beta }\int_0^{1-\frac{1}{k}}t^k\,dt \,\gtrsim \frac{1}{k(\log
k)^\beta }.$$
\par
Thus, we have seen that $\mu $ is a $\beta $-logarithmic
$1$-Carleson measure which satisfies
\begin{equation}\label{prec-moments}
\mu_n\,\asymp \,\frac{1}{n(\log n)^\beta }.\end{equation} Take $p\in
(1,\infty )$ and  $\alpha >\frac{1}{p}$ and set
$$a_n\,=\,\frac{1}{(n+1)\left (\log (n+2)\right)^\alpha },\quad
n\,=\,0,\,1,\,2,\,\dots ,$$ and
$$g(z)=\sum_{n=0}^\infty a_n\,z^n,\quad z\in
\mathbb D.$$ Notice that $\{ a_n\} \downarrow 0$ and that
 $\sum_{n=0}^\infty n^{p-1}\vert a_n\vert ^p\,<\,\infty $. Hence,
 $g\in B^p$.
\par We are going to prove that $\mathcal H_\mu (g)\,\not\in \,B^p$. This
implies that $\mathcal H_\mu (B^p)\,\not\subset B^p$, proving the
theorem.
\par We have
$\mathcal{H}_\mu (g)(z)=
\sum_{n=0}^{\infty}\left(\sum_{k=0}^{\infty}
\mu_{n+k}{a_k}\right)z^n$. Notice that $a_k\ge 0$ for all $k$ and
that the sequence of moments $\{ \mu _n\} $ is a decreasing sequence
of non-negative numbers. Then it follows that the sequence $\{
\sum_{k=0}^{\infty} \mu_{n+k}{a_k}\} _{n=0}^\infty $ of the Taylor
coefficients of $\mathcal H_\mu (g)$ is decreasing. Consequently, we
have that
\begin{equation}\label{gBp}\mathcal H_\mu (g)\,\in B^p\,\,\Leftrightarrow
\,\, \sum_{n=1}^\infty n^{p-1}\left \vert \sum_{k=0}^{\infty}
\mu_{n+k}{a_k}\right \vert ^p\,<\,\infty .\end{equation} Using the
definition of the sequence $\{ a_k\} $, (\ref{prec-moments}) and the
simple inequalities $\frac{k}{n+k}\ge \frac{1}{n+1}$ and ${\log
(n+k)}\le (\log n)(\log k)$ which hold whenever $k, n\,\ge 10$, say,
we obtain
\begin{align*} &\sum_{n=1}^\infty n^{p-1}\left \vert \sum_{k=0}^{\infty}
\mu_{n+k}{a_k}\right \vert ^p\,\ge \, \sum_{n=10}^\infty
n^{p-1}\left \vert \sum_{k=10}^{\infty} \mu_{n+k}{a_k}\right \vert
^p\\ \gtrsim & \sum_{n=10}^\infty n^{p-1}\left
(\sum_{k=10}^{\infty}\left [ \frac{1}{(n+k)\left (\log (n+k)\right
)^\beta}\,\frac{1}{k\left (\log k\right )^\alpha }\right ]\right )^p
\\ \gtrsim &
\sum_{n=10}^\infty \frac{1}{n(\log n)^{p\beta }}\left
(\sum_{k=10}^{\infty} \frac{1}{k^2\,\left (\log k\right )^{\alpha
+\beta }}\right )^p \,=\,\infty .
\end{align*}
Bearing in mind (\ref{gBp}), this implies that $H_\mu (g)\,\not\in
B^p$ as desired.
\end{Pf}
\par\bigskip
\begin{Pf}{\,\em{Theorem \ref{Bp-log-nec}.}}
Suppose that $1<p<\infty $ and $\gamma <1-\frac{1}{p}$. Let $\mu $
be a positive Borel measure on $[0,1)$ such that the operator
$\mathcal H_\mu $ is a bounded operator from $B^p$ into itself. Set
$\alpha \,=\,1-\gamma $,
$$a_k=\frac{1}{k(\log k)^\alpha },\quad k\ge 2,$$ and
$$f(z)=\sum_{k=2}^\infty a_k\,z^k,\quad z\in
\mathbb D.$$ Since $\alpha >\frac{1}{p}$, using
Theorem\,\@\ref{coef-dec-Besov} we see that $f\in B^p$. By our
assumption $H_\mu (f)\in B^P$, that is, $\Vert \mathcal H_\mu
(f)\Vert _{B^p}<\infty $. We have
$$\mathcal H_\mu (f)(z)\,=\,\sum_{n=0}^\infty \left
(\sum_{k=2}^\infty \mu_{n+k}a_k\right )z^n.$$ Since $a_k\ge 0$ for
all $k$ and $\{ \mu _n\} $ is a decreasing sequence of non-negative
numbers, it follows that the sequence $\{ \sum_{k=2}^\infty
\mu_{n+k}a_k\} _{n=0}^\infty $ is a decreasing sequence of
non-negative numbers. Then, using Theorem\,\@\ref{coef-dec-Besov} we
obtain
\begin{align*}
\Vert \mathcal H_\mu (f)\Vert _{B^p}^p\,\gtrsim \,&\sum_{n=1}^\infty
n^{p-1}\left (\sum_{k=2}^\infty \mu_{n+k}a_k\right )^p\\ \gtrsim \,&
\,\sum_{n=1}^\infty n^{p-1}\left (\sum_{k=2}^\infty \frac{1}{k(\log
k)^\alpha }\int_{[0,1)}x^{n+k}\,d\mu (x)\right )^p\\ \ge \,&
\sum_{n=1}^\infty n^{p-1}\left (\sum_{k=2}^\infty \frac{1}{k(\log
k)^\alpha }\int_{[t,1)}x^{n+k}\,d\mu (x)\right )^p
\\ \ge \,& \sum_{n=1}^\infty n^{p-1}\left (\sum_{k=2}^\infty \frac{t^{n+k}}{k(\log
k)^\alpha }\right )^p\,\mu \left ([t,1)\right )^p\\ = \,&
\sum_{n=1}^\infty n^{p-1}t^{np}\left (\sum_{k=2}^\infty
\frac{t^{k}}{k(\log k)^\alpha }\right )^p\,\mu \left ([t,1)\right
)^p,\quad\text{for all $t\in (0, 1)$}.
\end{align*}
Now, it is well known that $\sum_{k=2}^\infty \frac{t^{k}}{k(\log
k)^\alpha }\,\asymp \,\left (\log\frac{2}{1-t}\right )^{1-\alpha
}\,=\,\left (\log\frac{2}{1-t}\right )^{\gamma }$ (see
\cite[Vol.\,\@I, p.\,\@192]{Zy}). Then it follows that

\begin{align*}\Vert \mathcal H_\mu (f)\Vert _{B^p}^p\,\gtrsim &\,\left (\log\frac{2}{1-t}\right )^{\gamma p}
\left (\sum_{n=1}^\infty n^{p-1}t^{np}\right )\mu \left ([t,1)\right
)^p\\ \asymp &\, \left (\log\frac{2}{1-t}\right )^{\gamma p
}\frac{1}{(1-t)^p}\,\mu \left ([t,1)\right )^p.
 \end{align*}
Since $\Vert \mathcal H_\mu (f)\Vert _{B^p}<\infty $, this shows
that $\mu $ is a $\gamma $-logarithmic $1$-Carleson measure.
\end{Pf}
\par\bigskip
The following lemma will be used to prove
Theorem\,\@\ref{log-car-suf}. It is an adaptation of
\cite[Lemma\,\@7]{GaGiPeSis} to our setting. The proof is very
similar to that of the latter but we include it for the sake of
completeness.
\begin{lemma}\label{le:3} Let $p$, $\gamma $, and $\mu $ be as in
Theorem\,\@\ref{log-car-suf}. Then, there exists a constant $C=C(p,
\gamma , \mu )>0$ such that if $f\in B^p$, $g(z)=\sum_{k=0}^\infty
c_k z^k\in \hol(\D),$ and we set
$$h(z)=\sum_{k=0}^{\infty}c_k\left(\int_0^1t^{k+1}f(t)\,d\mu (t)\right)z^k,$$
then
\begin{equation*}
\|\Delta_n h\|_{H^p} \le C\left(\int_0^1 t^{2^{n-2}+1}|f(t)|\,d\mu
(t)\right)\| \Delta_ng\|_{H^p},\quad n\ge 3.
\end{equation*}
\end{lemma}

\begin{pf}
For  each $n=1,2, \dots$, define
$$
\Upsilon_n(s)=\int_0^1 t^{2^n s+1}f(t)\,d\mu (t),\quad s\ge 0.
$$
Clearly, $\Upsilon_n$ is a $C^\infty(0,\infty)$-function and
\begin{equation}\label{eq:up1}
\left| \Upsilon_n(s)\right|\le \int_0^1 t^{2^{n-2}+1}|f(t)|\,d\mu
(t),\quad s\ge \frac12.
\end{equation}

Furthermore, since $
 \sup_{0<x<1}\left(\log
\frac{1}{x}\right)^2 x^{1/2}=C(2)<\infty, $ we have
\begin{equation}
\begin{split}\label{eq:up2}
\left| \Upsilon''_n(s)\right| & \le
\int_0^1\left[\left(\log\frac{1}{t^{2^n}}\right)^2
t^{2^{n-1}}\right]\,t^{2^ns+1-2^{n-1}}|f(t)|\,d\mu (t)
\\  & \le C(2)\int_0^1 t^{2^ns+1-2^{n-1}}|f(t)|\,d\mu (t)
\le C(2) \int_0^1 t^{2^{n-2}+1}|f(t)|\,d\mu (t), \quad s\ge
\tfrac34.
\end{split}
\end{equation}

 Then, using (\ref{eq:up1}) and (\ref{eq:up2}), for each $n=1,2,\dots$, we can take
a function $\Phi_n\in C^\infty(\mathbb R)$ with
$\supp(\Phi_n)\in\left(\frac34, 4\right)$, and such that
$$
\Phi_n(s)= \Upsilon_n(s),\quad s\in\,[1,2],
$$
and
$$
A_{\Phi_n}=\max_{s\in\mathbb R}|\Phi_n(s)|+\max_{s\in\mathbb
R}|\Phi_n''(s)|\le C\int_0^1 t^{2^{n-2}+1}|f(t)|\,d\mu (t).
$$
Following the notation used in \cite[p.\,\@236]{GaGiPeSis}, we can
then write
\begin{align*}
\Delta_nh(z)&=\sum_{k=2^n}^{2^{n+1}-1}c_k\left(\int_0^1t^{k+1}f(t)\,d\mu
(t)\right)z^k \\ &=\sum_{k=2^n}^{2^{n+1}-1}c_k\Phi_n\left
(\frac{k}{2^n}\right )z^k=W_{2^n}^{\Phi_n}\ast\Delta_n g(z).
\end{align*}
So by using part (iii) of Theorem~B of \cite{GaGiPeSis}, we have
\begin{equation*}
\begin{split}
\| \Delta_n h\|_{H^p} &= \|W_{2^n}^{\Phi_{n}}\ast \Delta_n
g\|_{H^p}\leq C_pA_{\Phi_n}\| \Delta_n g\|_{H^p}\\ &\leq C
\left(\int_0^1 t^{2^{n-2}+1}|f(t)|\,d\mu (t)\right)\| \Delta_n
g\|_{H^p}.
\end{split}
\end{equation*}
\end{pf}

\par\bigskip
\begin{Pf}{\,\em{Theorem\,\@\ref{log-car-suf}.}}
By the closed graph theorem it suffices to show that $\mathcal H_\mu
(B^p)\subset B^p$.
\par Take $f\in B^p$. Since $\mu $ is a $\gamma $-logarithmic
$1$-Carleson measure, using Lemma\,\@\ref{HmuImuBp} we see that
\begin{equation*}\label{HI}\mathcal H_\mu (f)(z)\,=\,I_\mu (f)(z)\,=\,
\sum_{n=0}^\infty \left (\int_{[0,1)}t^nf(t)\,d\mu (t)\right
)z^n,\quad z\in \mathbb D.\end{equation*} Also, using
Corollary\,\@\ref{block-Besov}, we see that
\begin{equation}\label{Hmublock}
\mathcal H_\mu (f)\,\in \,B^p\,\,\,\Leftrightarrow \,\,\,
\sum_{n=1}^\infty 2^{-n(p-1)}\Vert \Delta _n\left (\mathcal H_\mu
(f)^\prime \right )\Vert _{H^p}^p\,<\,\infty .\end{equation} Now, we
have
\begin{equation*}\label{blockofHmuprime}
\Delta _n\left (\mathcal H_\mu (f)^\prime \right
)(z)\,=\,\sum_{k=2^n}^{2^{n+1}-1}(k+1)\left
(\int_{[0,1)}t^{k+1}f(t)\,d\mu (t)\right )z^k.\end{equation*} Using
Lemma\,\@\ref{le:3} we obtain that
\begin{equation*}\Vert \Delta _n\left (\mathcal H_\mu (f)^\prime \right
)\Vert _{H^p}\,\lesssim \,\left (\int_{[0,1)}t^{2^{n-2}+1}\vert
f(t)\vert \,d\mu (t)\right )\Vert \Delta
_nF\Vert_{H^p}\end{equation*} with $F(z)=\sum_{k=0}^\infty (k+1)z^k$
($z\in \mathbb D$). Now, we have that $M_p(r,F)\,=\,\og \left
(\frac{1}{(1-r)^{2-\frac{1}{p}}}\right )$ and then it follows that
$\Vert \Delta _nF\Vert_{H^p}\,=\,\og \left
(2^{n(2-\frac{1}{p})}\right )$ (see, e.\,\@g., \cite{MP-stu}). Using
this and the estimate $\vert f(t)\vert \lesssim \left (\log
\frac{2}{1-t}\right )^{1/p^\prime}$, we obtain
\begin{equation*}
\Vert \Delta _n\left (\mathcal H_\mu (f)^\prime \right )\Vert
_{H^p}\,\lesssim \,2^{n(2-\frac{1}{p})}\left
(\int_{[0,1)}t^{2^{n-2}+1} \left (\log \frac{2}{1-t}\right
)^{1/p^\prime}\,d\mu (t)\right ),\end{equation*} which using the
fact that $\mu $ is a $\gamma $-logarithmic $1$-Carleson measure and
Lemma\,\@\ref{gen-lem-moments} implies
\begin{equation*}\label{Hp-norm-Deltan-Hmuprime}
\Vert \Delta _n\left (\mathcal H_\mu (f)^\prime \right )\Vert
_{H^p}\,\lesssim \,2^{n(2-\frac{1}{p})}2^{-n}n^{\frac{1}{p^\prime
}-\gamma }\,=\,2^{n/p^\prime }n^{\frac{1}{p^\prime }-\gamma }.
\end{equation*}
This, together with the fact that $\gamma >1$, implies that
\begin{align*} &
\sum_{n=1}^\infty 2^{-n(p-1)}\Vert \Delta _n\left (\mathcal H_\mu
(f)^\prime \right )\Vert _{H^p}^p\,\lesssim \,\sum_{n=1}^\infty
2^{-n(p-1)}2^{np/p^\prime }n^{p(1-\gamma )-1}\\
& =\,\sum_{n=1}^\infty n^{p(1-\gamma )-1}\,<\,\infty .\end{align*}
Bearing in mind (\ref{Hmublock}), this shows that $\mathcal H_\mu
(f)\,\in \,B^p$ and finishes the proof.
\end{Pf}
\par\bigskip
\begin{center}{\bf Ackowledgements.}\end{center}
\par The authors wish to express their gratitude to the referee who read the paper very carefully and
made a good number of suggestions for improvement.
\par This research is supported by a grant from \lq\lq El Ministerio de
Econom\'{\i}a y Competitividad\rq\rq , Spain (MTM2014-52865-P) and
by a grant from la Junta de Andaluc\'{\i}a FQM-210. The second
author is also supported by a grant from \lq\lq El Ministerio de de
Educaci\'{o}n, Cultura y Deporte\rq\rq , Spain (FPU2013/01478).

\end{document}